\documentclass[12pt]{article}

\usepackage{amsmath}
\usepackage{amssymb}
\usepackage{braket} 
\usepackage{enumerate} 
\usepackage{amsthm}
\usepackage{graphicx}
\usepackage{psfrag}
\usepackage{framed}
\usepackage[T1]{fontenc}
\input{xy}  
\xyoption{all} 
\CompileMatrices 

\newtheorem{thm}{Theorem}[section]
\newtheorem{prop}[thm]{Proposition}
\newtheorem{lem}[thm]{Lemma}
\newtheorem{cor}[thm]{Corollary}

\numberwithin{equation}{section}

\theoremstyle{definition}
\newtheorem{dfn}[thm]{Definition}

\newtheorem{rem}[thm]{Remark}

\newcommand{\puk}[1]{\mathop{\mathrm{Puk}}\left(#1\right)}
\newcommand{\pukN}[2]{\mathop{\mathrm{Puk}}\nolimits_{#2}\left(#1\right)}
\newcommand{\type}[1]{\mathop{\mathrm{Type}}\left(#1\right)}
\newcommand{\Span}{\mathop{\mathrm{Span}}\nolimits}
\newcommand{\IIi}{$\textrm{II}_1\ $}
\newcommand{\tr}{\mathrm{tr}}
\newcommand{\ip}[2]{\left<{#1},{#2}\right>}

\newcommand{\ce}[2]{\mathbb E_{#1}\left(#2\right)}
\newcommand{\nm}[1]{\left\|{#1}\right\|}

\newcommand{\A}{\mathcal A}

\newcommand{\C}{\mathcal C}

\newcommand{\Ad}[1]{\textrm{Ad}(#1)}
\newcommand{\out}{\mathrm{Out}}

\newcommand{\VN}[1]{\mathcal L(#1)}
\newcommand{\vnotimes}{\overline{\otimes}}
\newcommand{\Ni}{\mathbb {N_\infty}}

\newcommand{\ISet}[2]{I(#1,#2)}

\title{Values of the Puk\'anszky invariant in McDuff factors}
\author{Stuart White\thanks{\texttt{s.white@maths.gla.ac.uk}}}
\date{}

\begin{document}

\maketitle

\begin{abstract}
In 1960 Puk\'anszky introduced an invariant associating to every masa in a separable \IIi factor a non-empty subset of $\mathbb N\cup\{\infty\}$.  This invariant examines the multiplicity structure of the von Neumann algebra generated by the left-right action of the masa.  In this paper it is shown that any non-empty subset of $\mathbb N\cup\{\infty\}$ arises as the Puk\'anszky invariant of some masa in a separable McDuff \IIi factor containing a masa with Puk\'anszky invariant $\{1\}$.  In particular the hyperfinite \IIi factor and all separable McDuff \IIi factors with a Cartan masa satisfy this hypothesis.  In a general separable McDuff \IIi factor we show that every subset of $\mathbb N\cup\{\infty\}$ containing $\infty$ is obtained as a Puk\'anszky invariant of some masa.
\end{abstract}

\section{Introduction}
In \cite{Pukanszky.Invariant} Puk\'anszky introduced an invariant for a maximal abelian self-adjoint subalgebra (masa) inside a separable \IIi factor, which he used to exhibit a countable infinite family of singular masas in the hyperfinite \IIi factor no pair of which are conjugate by an automorphism.  The invariant associates a non-empty subset of $\mathbb N\cup\{\infty\}$ to each masa $A$ in a separable \IIi factor $N$ as follows. Let $\A$ be the abelian von Neumann subalgebra of $\mathbb B(L^2(N))$ generated by $A$ and $JAJ$, where $J$ denotes the canonical involution operator on $L^2(N)$.  The orthogonal projection $e_A$ from $L^2(N)$ onto $L^2(A)$ lies in $\A$ and the algebra $\A'(1-e_A)$ is type $\mathrm{I}$ so decomposes as a direct sum of type $\mathrm{I}_{n}$-algebras.  The Puk\'anszky invariant of $A$ is the set of those $n\in\mathbb N\cup\{\infty\}$ appearing in this decomposition and is denoted $\puk{A}$.  See also \cite[Section 2]{Sinclair.Puk}.  

There has been recent interest in the range of values of the Puk\'anszky invariant in various \IIi factors.  Nesheyev and St{\o}rmer used ergodic constructions to show that any set containing $1$ arises as a Puk\'anszky invariant of a masa in the hyperfinite \IIi factor \cite[Corollary 3.3]{Stormer.Puk}.  Sinclair and Smith produced further subsets using group theoretic properties in \cite{Sinclair.Puk} and with Dykema in \cite{Sinclair.FreePuk}, which also examines free group factors. In the other direction Dykema has shown that $\sup\puk{A}=\infty$, whenever $A$ is a masa in a free group factor \cite{Dykema.AppFreeEntropy}.  

In this paper we show that every non-empty subset of $\mathbb N\cup\{\infty\}$ arises as the Puk\'anskzy invariant of a masa in the hyperfinite \IIi factor by means of an approximation argument.  More generally we obtain the same result in any separable McDuff \IIi factor containing a simple masa, that is one with Puk\'anskzy invariant $\{1\}$ (Corollary \ref{Main.SimpleMcDuff}).  These factors are the first for which the range of the Puk\'anskzy invariant has been fully determined.  Without assuming the presence of a simple masa we are able to show that every separable McDuff \IIi factor contains a masa with Puk\'anskzy invariant $\{\infty\}$ and hence we obtain every subset of $\mathbb N\cup\{\infty\}$ containing $\infty$ as a Puk\'anszky invariant of some masa in these factors (Theorem \ref{Main.GeneralMcDuff}).  In particular, there are uncountably many singular masas in any separably McDuff factor, no pair of which is conjugate by an automorphism of the factor.

Section \ref{Values} contains a construction for producing masas in McDuff \IIi factors. Given a McDuff \IIi factor $N_0$ we shall repeatedly tensor on copies of the hyperfinite \IIi factor --- this gives us a chain $(N_s)_{s=0}^\infty$ of \IIi factors whose direct limit $N$ is isomorphic to $N_0$.  We shall produce a masa $A$ in $N$ by giving an approximating sequence of masas $A_s$ in each $N_s$ such that $A_s\subset A_{s+1}$ and defining $A=(\bigcup_{s=0}^\infty A_s)''$.  This idea has its origin in \cite{Tauer.Masa} working in the hyperfinite \IIi factor arising as the infinite tensor produce of finite matrix algebras, although using finite matrix algebras can only yield masas with Puk\'anskzy invariant $\{1\}$, \cite[Theorem 4.1]{Saw.Tauer}.

In the remainder of the introduction we outline the construction of a masa with Puk\'anszky invariant $\{2,3\}$.  Initially we shall produce a masa $A_1$ in $N_1$ such that the multiplicity structure of $\A_1$ (the algebra generated by the left-right action of $A_1$ on $L^2(N_1)$) is represented by Figure \ref{Intro.Fig1}.  By this we mean that $e$ is a projection of trace $1/2$ in $A$ and that $\A_1'eJeJ$ and $\A_1'e^\perp Je^\perp J$ are both type $\mathrm{I}_1$, while $\A_1'eJe^\perp J$ and $\A_1'e^\perp JeJ$ are type $\mathrm{I}_2$.
\begin{figure}[h!]
\begin{center}
\begin{picture}(100,90)
\put(20,0){\line(1,0){80}}
\put(20,40){\line(1,0){80}}
\put(20,80){\line(1,0){80}}
\put(20,0){\line(0,1){80}}
\put(60,0){\line(0,1){80}}
\put(100,0){\line(0,1){80}}
\put(37,17){2}
\put(77,57){2}
\put(37,57){1}
\put(77,17){1}
\put(37,83){$e$}
\put(77,83){$e^\perp$}
\put(0,57){$e$}
\put(0,17){$e^\perp$}
\end{picture}
\caption{Symbolic description of the multiplicity structure of $\A_1$.}
\label{Intro.Fig1}
\end{center}
\end{figure}

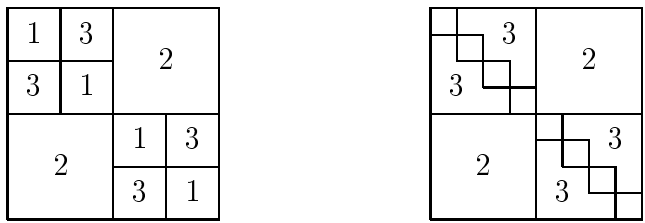
\begin{figure}[h!]
\begin{center}
\begin{picture}(240,80)
\put(0,0){\line(1,0){80}}
\put(0,40){\line(1,0){80}}
\put(0,80){\line(1,0){80}}

\put(0,0){\line(0,1){80}}
\put(40,0){\line(0,1){80}}
\put(80,0){\line(0,1){80}}
\put(17,17){2}
\put(57,57){2}
\put(60,0){\line(0,1){40}}
\put(0,60){\line(1,0){40}}
\put(40,20){\line(1,0){40}}
\put(20,40){\line(0,1){40}}

\put(7,47){3}
\put(27,67){3}
\put(7,67){1}
\put(27,47){1}
\put(47,7){3}
\put(67,27){3}
\put(47,27){1}
\put(67,7){1}

\put(160,0){\line(1,0){80}}
\put(160,40){\line(1,0){80}}
\put(160,80){\line(1,0){80}}
\put(160,0){\line(0,1){80}}
\put(200,0){\line(0,1){80}}
\put(240,0){\line(0,1){80}}
\put(177,17){2}
\put(217,57){2}
\put(230,0){\line(0,1){20}}
\put(220,10){\line(0,1){20}}
\put(210,20){\line(0,1){20}}
\put(200,30){\line(0,1){10}}

\put(160,70){\line(1,0){20}}
\put(170,60){\line(1,0){20}}
\put(180,50){\line(1,0){20}}
\put(190,40){\line(1,0){10}}

\put(200,30){\line(1,0){20}}
\put(210,20){\line(1,0){20}}
\put(220,10){\line(1,0){20}}
\put(230,0){\line(1,0){10}}

\put(190,40){\line(0,1){20}}
\put(180,50){\line(0,1){20}}
\put(170,60){\line(0,1){20}}
\put(160,70){\line(0,1){10}}

\put(167,47){3}
\put(187,67){3}
\put(207,7){3}
\put(227,27){3}
\end{picture}
\caption{The multiplicity structures of $\A_2$ and $\A_3$.}
\label{Intro.Fig2}
\end{center}
\end{figure}

At the second stage we subdivide $e$ and $e^\perp$ to obtain four projections in $A_2$ and arrange for the multiplicity structure of $\A_2$ to be represented by the left diagram in Figure \ref{Intro.Fig2}.  We then cut each of these projections in half again and ensure that the multiplicity structure of $\A_3$ is represented by the second diagram in Figure \ref{Intro.Fig2}, where $1$'s appear down the diagonal.  It is important to do this in such a way that a limiting argument can be used to obtain the multiplicity structure of $\A=(A\cup JAJ)''$.  If this is done successfully, then the multiplicity structure of $\A$ will be represented by Figure \ref{Intro.Fig3}, where the diagonal line has multiplicity $1$.
If we further ensure that the projections used to cut down the masas $A_r$ in this construction generate $A$, then the diagonal line in Figure \ref{Intro.Fig3} corresponds to the projection $e_A$ with range $L^2(A)$ and this is the projection explicitly removed in the definition of $\puk{A}$.  The resulting masa $A$ will then have Puk\'anszky invariant $\{2,3\}$ as required.  

\begin{figure}[h!]
\begin{center}
\begin{picture}(80,80)
\put(0,0){\line(1,0){80}}
\put(0,40){\line(1,0){80}}
\put(0,80){\line(1,0){80}}
\put(0,0){\line(0,1){80}}
\put(40,0){\line(0,1){80}}
\put(80,0){\line(0,1){80}}
\put(0,80){\line(1,-1){80}}
\put(17,17){2}
\put(57,57){2}
\put(27,67){3}
\put(7,47){3}
\put(47,7){3}
\put(67,27){3}
\end{picture}
\caption{The multiplicity structure of $\A$.}
\label{Intro.Fig3}
\end{center}
\end{figure}

To get from Figure \ref{Intro.Fig1} to the left diagram in Figure \ref{Intro.Fig2} in a compatible way, we `tensor on' the diagram in Figure \ref{Intro.Fig4}.  This is done by producing masas $D_1,D_2,D_3,D_4$ in the hyperfinite \IIi factor $R$ such that $(D_i\cup JD_jJ)'$ is type $\mathrm{I}_1$ unless $i,j$ is the unordered pair $\{1,2\}$ or $\{3,4\}$.  In these cases $(D_i\cup JD_jJ)'$ is type $\mathrm{I}_3$.  Given projections $e_1,e_2,e_3,e_4$ in $A_1$ with $e=e_1+e_2$ and $e^\perp=e_3+e_4$ and $\tr(e_i)=1/4$ for each $i$ we shall define $A_2$ in $N_2=N_1\ \vnotimes\ R$ by
$$
A_2=\bigoplus_{i=1}^4A_1e_i\ \vnotimes D_i.
$$
In this way $\A_2$ has the required multiplicity structure.
\begin{figure}[h!]
\begin{center}
\begin{picture}(80,80)
\put(0,0){\line(1,0){80}}
\put(0,20){\line(1,0){80}}
\put(0,40){\line(1,0){80}}
\put(0,60){\line(1,0){80}}
\put(0,80){\line(1,0){80}}

\put(0,0){\line(0,1){80}}
\put(20,0){\line(0,1){80}}
\put(40,0){\line(0,1){80}}
\put(60,0){\line(0,1){80}}
\put(80,0){\line(0,1){80}}
\put(7,7){1}
\put(27,7){1}
\put(7,27){1}
\put(27,27){1}
\put(47,47){1}
\put(67,67){1}
\put(47,67){1}
\put(67,47){1}
\put(7,47){3}
\put(27,67){3}
\put(7,67){1}
\put(27,47){1}
\put(47,7){3}
\put(67,27){3}
\put(47,27){1}
\put(67,7){1}
\end{picture}
\caption{Mixed Puk\'anszky invariant structure of the masas $D_1,D_2,D_3,D_4$.}\label{Intro.Fig4}
\end{center}
\end{figure}

In sections \ref{Prelim} and \ref{SumCartan} we develop the concept of mixed Puk\'anszky invariants of pairs of masas to handle the families $(D_i)$, which we will repeatedly adjoin.  The main result is Theorem \ref{SumCartan.CountableFamily}, which ensures that the family $D_1,D_2,D_3,D_4$ above, and other families in this style can indeed be found.  In section \ref{Values} we give the details of the inductive construction and in section \ref{MainLem} we compute the Puk\'anszky invariant of the resulting masa.  We end in section \ref{Main} by collecting together the main results.

\section{Mixed Puk\'anszky Invariants}\label{Prelim}
In this paper all \IIi factors will be separable.  In this way we only need one infinite cardinal denoted $\infty$.  We shall write $\Ni$ for the set $\mathbb N\cup\{\infty\}$ henceforth.
\begin{dfn}
Given a type $\mathrm{I}$ von Neumann algebra $M$ we shall write $\type{M}$ for the set of those $m\in\Ni$ such that $M$ has a non-zero component of type $\mathrm{I}_m$.
\end{dfn}
Given a \IIi factor $N$, write $\tr$ for the unique faithful trace on $N$ with $\tr(1)=1$.  For $x\in N$, let $\nm{x}_2=\tr(x^*x)^{1/2}$, a pre-Hilbert space norm on $N$.  The completion of $N$ in this norm is denoted $L^2(N)$.  Define a conjugate linear isometry $J$ from $L^2(N)$ into itself by extending $x\mapsto x^*$ by continuity from $N$.

\begin{dfn}\label{Prelim.DefMix}
Given two masas $A$ and $B$ in a \IIi factor $N$ define the \emph{mixed Puk\'anszky invariant} of $A$ and $B$ to be the set $\type{(A\cup JBJ)'}$, where the commutant is taken in $\mathbb B(L^2(N))$.  We denote this set $\puk{A,B}$ or $\pukN{A,B}{N}$ when it is necessary. Note that $\puk{A,A}=\puk{A}\cup\{1\}$ for any masa $A$, the extra $1$ arising as the Jones projection $e_A$ is not removed in the definition of $\puk{A,A}$.  
\end{dfn}

It is immediate that $\puk{A,B}$ is a conjugacy invariant of a pair of masas $(A,B)$ in a \IIi factor, i.e. that if $\theta$ is an automorphism of $N$ we have $\puk{A,B}=\puk{\theta(A),\theta(B)}$.  If we only apply $\theta$ to one masa in the pair then we may get different mixed invariants.  For an inner automorphism this is not the case.
\begin{prop}\label{Prelim.Inner}
Let $A$ and $B$ be masas in a \IIi factor $N$.  For any unitaries $u,v\in N$ we have
$$
\puk{uAu^*,vBv^*}=\puk{A,B}.
$$
\end{prop}
\begin{proof}
Consider the automorphism $\Theta=\Ad{uJvJ}$ of $B(L^2(N))$, which has $\Theta(A)=uAu^*$ and $\Theta(JBJ)=JvBv^*J$.  Therefore $(A\cup JBJ)'$ and $(uAu^*\cup J(vBv^*)J)'$ are isomorphic, so have the same type decomposition.\qed
\end{proof}

The Puk\'anskzy invariant is well behaved with respect to tensor products \cite[Lemma 2.1]{Sinclair.Puk}.  So too is the mixed Puk\'anszky invariant.  Given $E,F\subset\Ni$ write $E\cdot F=\Set{mn|m\in E,n\in F}$, where by convention $n\infty=\infty n=\infty$ for any $n\in\Ni$.  
\begin{lem}\label{Prelim.MixTensor}
Let $(N_i)_{i\in I}$ be a countable family of finite factors.  Suppose that we have masas $A_i$ and $B_i$ in $N_i$ for each $i\in I$.  Let $N$ be the finite factor obtained as the infinite von Neumann tensor product of the $N_i$ with respect to the product trace and let $A$ and $B$ be the infinite tensor products of the $A_i$ and $B_i$ respectively.  Then $A$ and $B$ are masas in $N$. When $I$ is finite,
$$
\pukN{A,B}{N}=\prod_{i\in I}\pukN{A_i,B_i}{N_i}.
$$
If $I$ is infinite, and each $\pukN{A_i,B_i}{N_i}=\{n_i\}$ for some $n_i\in\Ni$, then $\pukN{A,B}{N}=\{n\}$, where $n=\prod_In_i$, when all but finitely many $n_i=1$, and $n=\infty$ otherwise.
\end{lem}
\begin{proof}
That $A$ and $B$ are masas follows from Tomita's commutation theorem, see \cite[Theorem 11.2.16]{KR.2}.  Suppose first that $I$ is finite. For each $i\in I$, let $(p_{i,n})_{n\in\Ni}$ be the decomposition of the identity projection into projections in $(A_i\cup JB_iJ)''\subset B(L^2(N_i))$ such that $(A_i\cup JB_iJ)'p_{i,n}$ is type ${\textrm{I}}_n$ for each $n\in\Ni$ (some of these projections may be zero). Then given any family $(n_i)_i$ in $\Ni$, $p=\bigotimes_{i\in I}p_{i,n_i}$ is a central projection in $(A\cup JBJ)'$ and $(A\cup JBJ)'p$ is type $\mathrm{I}_m$ where $m=\prod_{i\in I}n_i$.  All these projections are mutually orthogonal with sum $1$. Therefore $\pukN{A,B}{N}$ consists of those $m$ such that $p\neq 0$ and this occurs if and only if all the corresponding $p_{i,n_i}$ appearing in the tensor product are non-zero.  These are precisely the $m$ in $\prod_{i\in I}\pukN{A_i,B_i}{N_i}$.

Suppose $I$ is infinite and each $\pukN{A_i,B_i}{N_i}=\{n_i\}$, for some $n_i\in\Ni$.  Let $\A_i=(A_i\cup JB_iJ)''\subset \mathbb B(L^2(N_i))$ and $\A_i'$ it's commutant of $\A_i$ in $\mathbb B(L^2(N_i))$.  Let $\A=(A\cup JBJ)''$ in $\mathbb B(L^2(N))$ and $\A'$ the commutant of $\A$ in this algebra.  The Tomita commutation theorem gives 
$$
\A'=\overline{\bigotimes}\A_i'\subseteq\overline{\bigotimes}\mathbb B(L^2(N_i))\cong\mathbb B(L^2(N)).
$$
Since each $\A_i'\cong\A_i\vnotimes\mathbb M_{n_i}$, where $\mathbb M_{n_i}$ is the $n_i\times n_i$ matrices (or $\mathbb B(H)$ for some separable infinite dimensional Hilbert space when $n_i=\infty$).  Thus 
$$
\A'\cong\left(\overline{\bigotimes}A_i\right)\vnotimes\left(\overline{\bigotimes}\mathbb M_{n_i}\right)\cong A\vnotimes\mathbb M_n,
$$
so $\A'$ is homogenous of type ${\mathrm{I}}_n$.\qed
\end{proof}

Given two masas $A$ and $B$ in a \IIi factor $N$ we can form the algebra $M_2(N)$ of $2\times 2$ matrices over $N$.  We can construct a masa in $M_2(N)$ 
$$
\begin{pmatrix}A&0\\0&B\end{pmatrix}=\Set{\begin{pmatrix}a&0\\0&b\end{pmatrix}|a\in A,b\in B},
$$
which we denote $A\oplus B$ --- the direct sum of $A$ and $B$.   In \cite{Sinclair.Puk} it is noted that if $B$ is a unitary conjugate of $A$, then the Puk\'anszky invariant of $A\oplus B$ can be determined from that of $A$ (and hence $B$).  Indeed we have
$$
\puk{A\oplus uAu^*}=\puk{A}\cup\{1\},
$$
whenever $u$ is a unitary in $N$.  The initial motivation for the introduction of the mixed Puk\'anszky invariant was to aid in the study of the Puk\'anszky invariant of these direct sums since
$$
\puk{A\oplus B}=\puk{A}\cup\puk{B}\cup\puk{A,B},
$$
whenever $A$ and $B$ are masas in a \IIi factor $N$.  As we shall subsequently see, the Puk\'anskzy invariant behaves badly with respect to the direct sum construction.  In the next section we shall give Cartan masas $A$ and $B$ in the hyperfinite \IIi factor such that $\puk{A\oplus B}=\{1,n\}$ for any $n\in\Ni$, and given non-empty sets $E,F,G\subset\Ni$ we shall construct, in Theorem \ref{Main.DirectSum}, masas $A$ and $B$ in the hyperfinite \IIi factor such that $\puk{A}=E$, $\puk{B}=F$ and $\puk{A,B}=G$. Hence it is not possible to make a more general statement about the Puk\'anszky invariant of a direct sum than 
$$
\puk{A\oplus B}\supset\puk{A}\cup\puk{B}.
$$

\section{Mixed invariants of Cartan masas in $R$}\label{SumCartan}
In this section we shall construct large families of Cartan masas in the hyperfinite \IIi factor, each masa will have Puk\'anszky invariant $\{1\}$ by virtue of being Cartan \cite[Section 3]{Popa.NotesCartan}.  Our objective will be to control the mixed Puk\'anskzy invariant of any two elements from the family.  We start by constructing a family of three Cartan masas in the hyperfinite \IIi factor and then use Lemma \ref{Prelim.MixTensor} to produce the desired result.
\begin{lem}\label{SumCartan.ABC}
For each $n\in\Ni$ there exists Cartan masas $A,B,C$ in the hyperfinite \IIi factor such that $\puk{A,B}=\{n\}$ while $\puk{A,C}=\puk{B,C}=\{1\}$.
\end{lem}
We shall first establish Lemma \ref{SumCartan.ABC} when $n$ is finite.  The lemma is immediate for $n=1$, take $A=B=C$ to be any Cartan masa in the hyperfinite \IIi factor.  Let $n\geq 2$ be a fixed integer until further notice.  Since any two Cartan masas in the hyperfinite \IIi factor are conjugate by an automorphism \cite{Connes.CFW}, we shall fix a Cartan masa $A$ arising as the diagonals in an infinite tensor product and then construct $B=\theta(A)$ and $C=\phi(A)$ by exhibiting appropriate automorphisms $\theta$ and $\phi$ of $R$.  Let $M$ denote the $n\times n$ matrices and $D_0$ denote the diagonal $n\times n$ matrices, a masa in $M$.  Write $(e_i)_{i=0}^{n-1}$ for the minimal projections of $D_0$ so $e_i$ has $1$ in the $(i,i)$th entry and $0$ elsewhere.  Let 
$$
w=\begin{pmatrix}0&1&0&\dots&0\\0&0&1&\ddots&0\\ \vdots&\ddots&\ddots&\ddots&\vdots\\ 0&\ddots&\ddots&\ddots&1\\1&0&\dots&0&0\end{pmatrix}
$$
a unitary in $M$, which, in its action by conjugation, cyclically permutes the minimal projections of $D_0$.  That is $we_iw^*=e_{i-1}$ with the subtraction taken $\mod n$.  The abelian algebra generated by $w$ is a masa $D_1$ in $M$, which is orthogonal to $D_0$ \cite[Section 3]{Popa.Orth}.  Write $(f_i)_{i=0}^{n-1}$ for the minimal projections of $D_1$.  Define
\begin{equation}\label{SumCartan.Defv}
v=\sum_{i=0}^{n-1}w^i\otimes f_i
\end{equation}
a unitary in $D_1\otimes D_1\subset M\otimes M$.

We shall produce $A,B$ and $C$ in the hyperfinite \IIi factor $R$ realised as $(\bigotimes_{r=1}^\infty M)''$.  Let $A=(\bigotimes_{r=1}^\infty D_0)''$.  For each $r$ consider the unitary $u_r=1^{\otimes(r-1)}\otimes v$, which lies in $M^{\otimes(r+1)}\subset R$.  All of these unitaries commute (as they lie in the masa $(\bigotimes_{r=1}^\infty D_1)''$ in $R$) and satisfy $u_r^n=1$.  We are able to define automorphisms 
$$
\theta=\lim_{r\rightarrow\infty}\Ad{u_1u_2\dots u_r},\quad \phi=\lim_{r\rightarrow\infty}\Ad{u_1u_3u_5\dots u_{2r+1}}
$$
of $R$ with the limit taken pointwise in $\nm{.}_2$. Convergence follows, since for $x\in M^{\otimes r}$ we have $u_sxu_s^*=x$ whenever $s>r$ and such $x$ are $\nm{.}_2$-dense in $R$.  In this way $\theta$ and $\phi$ define $^*$-isomorphisms of $R$ into $R$.  As $\theta^n=I$ and $\phi^n=I$ (since the $u_r$s commute and each $u_r^n=1$), we see that $\theta$ and $\phi$ are onto and so automorphisms of $R$.  Define Cartan masas $B=\theta(A)$ and $C=\phi(A)$ in $R$.  The calculations of $\puk{A,C}$ and $\puk{B,C}$ are straightforward.

\begin{lem}\label{SumCartan.ABC.Easy}
With the notation above, we have $\puk{A,C}=\puk{B,C}=\{1\}$.
\end{lem}
\begin{proof}
We re-bracket the infinite tensor product defining $R$ as
$$
R=(M\otimes M)\vnotimes(M\otimes M)\vnotimes\dots
$$
so that $R$ is the infinite tensor product of copies of $M\otimes M$.  Since $u_{2r+1}$ lies in $1^{\otimes 2r}\otimes(M\otimes M)$ we see that $\phi$ factorises as $\prod_{s=1}^\infty\Ad{v}$ with respect to this decomposition.  Lemma \ref{Prelim.MixTensor} then tells us that $\puk{A,C}$ is the set product of infinitely many copies of $\pukN{D_0\otimes D_0,v(D_0\otimes D_0)v^*}{M\otimes M}$.  Since $D_0\otimes D_0$ and $v(D_0\otimes D_0)v^*$ are masas in $M_0\otimes M_0$ a simple dimension check ensures that $\pukN{D_0\otimes D_0,v(D_0\otimes D_0)v^*}{M\otimes M}=\{1\}$ and hence $\puk{A,C}=\{1\}$.

Observe that $\puk{B,C}=\puk{\theta(A),\phi(A)}=\puk{\phi^{-1}\theta(A),A}$.  As all the $u_r$ commute, we have 
$$
\phi^{-1}\circ\theta=\lim_{r\rightarrow\infty}\Ad{u_2u_4\dots u_{2r}}
$$
with pointwise $\nm{.}_2$ convergence.  This time we re-bracket the tensor product defining $R$ as
$$
R=M\vnotimes(M\otimes M)\vnotimes(M\otimes M)\vnotimes\dots,
$$
and since $u_{2r}=1^{\otimes 2r-1}\otimes v\in 1\otimes1^{\otimes{2(r-1)}}\otimes(M\otimes M)$, we obtain $\puk{B,C}=\{1\}$ in the same way.\qed
\end{proof}

The key tool in establishing that $\puk{A,B}=\{n\}$ is the following calculation, which we shall use to produce $n$ equivalent abelian projections for the commutant of the left-right action.
\begin{lem}\label{SumCartan.ABC.Keyclaim}
Use the notation preceding Lemma \ref{SumCartan.ABC.Easy}.  For $r=0,1,\dots,n-1$ let $\xi_r$ denote $f_r$ taken in the first copy of $M$ in the tensor product making up $R$, thought of as a vector in $L^2(R)$. For any $m\geq 0$, $i_1,i_2,\dots,j_m,j_1,j_2,\dots,j_m=0,1,\dots,n-1$ and $r,s=0,1,\dots,n-1$ we have
\begin{equation}\label{SumCartan.ABC.Keyclaim.1}
{\ip{(e_{i_1}\otimes\dots\otimes e_{i_m})\xi_r\theta(e_{j_1}\otimes\dots\otimes e_{j_m})}{\xi_s}}_{L^2(R)}=\delta_{r,s}n^{-(2m+1)}.
\end{equation}
\end{lem}
\begin{proof}
We proceed by induction.  When $m=0$, (\ref{SumCartan.ABC.Keyclaim.1}) reduces to $\ip{\xi_r}{\xi_s}=\delta_{r,s}n^{-1}$, which follows as $\ip{\xi_r}{\xi_s}=\tr(f_rf_s^*)$ and $(f_r)_{r=0}^{n-1}$ are the minimal projections of a masa in the $n\times n$ matrices.

For $m>0$ observe that $\theta(e_{j_1}\otimes\dots\otimes e_{j_m})=u_1\dots u_m(e_{j_1}\otimes\dots\otimes e_{j_m})u_m^*\dots u_1^*$. With the subtraction in the subscript taken $\mathrm{mod}\ n$, we have
$$
u_m(e_{j_1}\otimes\dots\otimes e_{j_m})u_m^*=e_{j_1}\otimes\dots\otimes e_{j_{m-1}}\otimes\Big(\sum_{k=0}^{n-1}e_{j_m-k}\otimes f_k\Big)
$$
from (\ref{SumCartan.Defv}) and $we_{j_m}w^*=e_{j_m-1}$.  Therefore
\begin{align}
&{\ip{(e_{i_1}\otimes\dots\otimes e_{i_m})\xi_r\theta(e_{j_1}\otimes\dots\otimes e_{j_m})}{\xi_s}}\nonumber\\
=&{\ip{(e_{i_1}\otimes\dots\otimes e_{i_m})\xi_ru_1\dots u_{m-1}\Big(e_{j_1}\otimes\dots\otimes e_{j_m-1}\otimes\sum_{k=0}^{n-1}e_{j_m-k}\otimes f_k\Big)u_{m-1}^*\dots u_1^*}{\xi_s}}\nonumber\\
=&\tr\Bigg(\sum_{k=0}^{n-1}\Big((e_{i_1}\otimes\dots\otimes e_{i_m})f_ru_1\dots u_{m-1}(e_{j_1}\otimes\dots \otimes e_{j_{m-1}}\otimes e_{j_m-k})u_{m-1}^*\dots u_1^*f_s^*\Big)\otimes f_k\Bigg)\nonumber\\
=&n^{-1}\tr\Bigg((e_{i_1}\otimes\dots\otimes e_{i_m})f_ru_1\dots u_{m-1}(e_{j_1}\otimes\dots \otimes e_{j_{m-1}}\otimes \sum_{k=0}^{n-1}e_{j_m-k})u_{m-1}^*\dots u_1^*f_s^*\Bigg)\nonumber\\
=&n^{-1}\tr\Bigg((e_{i_1}\otimes\dots\otimes e_{i_m})f_r\theta(e_{j_1}\otimes\dots \otimes e_{j_{m-1}}\otimes 1)f_s^*\Bigg)\label{SumCartan.ABC.KeyClaim.3}
\end{align}
as the $f_k$ in the third line is the only object appearing in the $(m+1)$-tensor position and $\tr$ is a product trace. This produces the factor $n^{-1}=\tr(f_k)$.  We obtain (\ref{SumCartan.ABC.KeyClaim.3}) as $e_{j_1}\otimes\dots\otimes e_{j_{m-1}}\otimes 1$ lies in $M^{\otimes(m-1)}$ so $\theta(e_{j_1}\otimes\dots\otimes e_{j_{m-1}}\otimes 1)=u_1\dots u_{m-1}(e_{j_1}\otimes\dots\otimes e_{j_{m-1}}\otimes 1)u_{m-1}^*\dots u_1^*$.

Now $\theta(f_r)=f_r$ for all $r$ (since each $u_m$ commutes with $f_r$) and $\theta$ is trace preserving.  In this way we obtain
\begin{align*}
&{\ip{(e_{i_1}\otimes\dots\otimes e_{i_m})\xi_r\theta(e_{j_1}\otimes\dots\otimes e_{j_m})}{\xi_s}}\\
=&n^{-1}\tr\Bigg(\theta^{-1}(e_{i_1}\otimes\dots\otimes e_{i_m})f_r(e_{j_1}\otimes\dots \otimes e_{j_{m-1}}\otimes 1)f_s^*\Bigg).
\end{align*}
We now apply the same argument again giving us
\begin{align*}
&{\ip{(e_{i_1}\otimes\dots\otimes e_{i_m})\xi_r\theta(e_{j_1}\otimes\dots\otimes e_{j_m})}{\xi_s}}\\
=&n^{-2}\tr\Bigg(\theta^{-1}(e_{i_1}\otimes\dots\otimes e_{i_{m-1}}\otimes 1)f_r(e_{j_1}\otimes\dots \otimes e_{j_{m-1}}\otimes 1)f_s^*\Bigg)\\
=&n^{-2}\tr\Bigg((e_{i_1}\otimes\dots\otimes e_{i_{m-1}})f_r\theta(e_{j_1}\otimes\dots \otimes e_{j_{m-1}})f_s^*\Bigg)\\
=&n^{-2}{\ip{(e_{i_1}\otimes\dots\otimes e_{i_{m-1}})\xi_r\theta(e_{j_1}\otimes\dots\otimes e_{j_{m-1}})}{\xi_s}}.
\end{align*}
The lemma now follows by induction.\qed
\end{proof}

We can now complete the proof of Lemma \ref{SumCartan.ABC}.
\begin{proof}{Proof of Lemma \ref{SumCartan.ABC}.}
We continue to let $n\geq 2$ be a fixed integer and let $A$ and $B$ be the masas introduced before Lemma \ref{SumCartan.ABC.Easy}.   Let $\C$ be the abelian algebra $(A\cup JBJ)''$ in $\mathbb B(L^2(R))$.  We continue to write $\xi_r$ for $f_r$ (in the first tensor position) thought of as a vector in $L^2(R)$.  For each $r$, let $P_r$ be the orthogonal projection in $\mathbb B(L^2(R))$ onto $\overline{\C\xi_r}$, an abelian projection in $\C'$.

Since elements $(e_{i_1}\otimes\dots\otimes e_{i_m})f_r\theta(e_{j_1}\otimes\dots\otimes e_{j_m})$, where $m\geq 0$ and $i_1,\dots,i_m,j_1,\dots,j_m=0,1,\dots,n-1$, have dense linear span in $\overline{\C\xi_r}$, Lemma \ref{SumCartan.ABC.Keyclaim} implies that $P_r$ is orthogonal to $P_s$ when $r\neq s$.  Furthermore, for each $m$, the elements 
$$
(e_{i_1}\otimes\dots\otimes e_{i_m})f_r\theta(e_{j_1}\otimes\dots\otimes e_{j_{m-1}}\otimes 1)
$$
indexed by $i_1,\dots,i_m,j_1,\dots,j_{m-1},r=0,1,\dots,n-1$ are $n^{2m}$ pairwise orthogonal non-zero elements of $M^{\otimes m}$, the $n^m\times n^m$ matrices.  Therefore, $M^{\otimes m}$ is contained in the range of $P_0+P_1+\dots+P_{n-1}$ for each $m$ so that $\sum_{r=0}^{n-1}P_r=1$.

It remains to show that all the $P_r$ are equivalent in $\C'$, from which it follows that $\C'$ is homogeneous of type $\mathrm{I}_n$.  Given $r\neq s$ we must define a partial isometry $v_{r,s}\in\C'$ with $v_{r,s}v_{r,s}^*=P_s$ and $v_{r,s}^*v_{r,s}=P_r$.  Lemma \ref{SumCartan.ABC.Keyclaim} allows us to define $v_{r,s}$ by extending the map $\xi_r\mapsto\xi_s$ by $(A,B)$-modularity.  More precisely define linear maps
$$
v_{r,s}^{(m)}:\Span(D_0^{\otimes m}f_r\theta(D_0^{\otimes m}))\rightarrow\Span(D_0^{\otimes m}f_s\theta(D_0^{\otimes m}))
$$
by extending
$$
v_{r,s}^{(m)}\Big((e_{i_1}\otimes\dots\otimes e_{i_m})f_r\theta(e_{j_1}\otimes\dots\otimes e_{j_m})\Big)=(e_{i_1}\otimes\dots\otimes e_{i_m})f_s\theta(e_{j_1}\otimes\dots\otimes e_{j_m})
$$
by linearity.  Lemma \ref{SumCartan.ABC.Keyclaim} shows that these maps preserve $\nm{.}_2$ and that $v_{r,s}^{(m+1)}$ extends $v_{r,s}^{(m)}$.  Let $v_{r,s}$ be the closure of the union of the $v_{r,s}^{(m)}$.  This is patently a partial isometry in $\C'$ with domain projection $P_r$ and range projection $P_s$.  Hence $\puk{A,B}=\{n\}$ and combining this with Lemma \ref{SumCartan.ABC.Easy} establishes Lemma \ref{SumCartan.ABC} when $n$ is finite.

When the $n$ of Lemma \ref{SumCartan.ABC} is $\infty$ we take a tensor product.  More precisely find Cartan masas $A_0,B_0,C_0$ in the hyperfinite \IIi factor $R_0$ such that $\puk{A_0,B_0}=\{2\}$ and $\puk{A_0,C_0}=\puk{B_0,C_0}=\{1\}$.  Now form the hyperfinite \IIi factor $R$ by taking the infinite tensor product of copies of $R_0$.  The Cartan masas $A$, $B$ and $C$ in $R$ obtained from the infinite tensor product of copies of $A_0$, $B_0$ and $C_0$ have $\puk{A,B}=\{\infty\}$, and $\puk{A,C}=\puk{B,C}=\{1\}$ by Lemma \ref{Prelim.MixTensor}.\qed
\end{proof}

\begin{rem}
By fixing a Cartan masa $D$ in a \IIi factor $N$ we could consider the map $\theta\mapsto\puk{D,\theta(D)}$, which (by Proposition \ref{Prelim.Inner}) induces a map on $\out{N}$.  This map is not necessarily constant on outer conjugacy classes, as the automorphisms $\theta$ and $\phi$ of the hyperfinite \IIi factor above have outer order $n$ and obstruction to lifting $1$ so are outer conjugate by \cite{Connes.OuterConjugacy}.
\end{rem}

Let us now give the main result of this section.
\begin{thm}\label{SumCartan.CountableFamily}
Let $I$ be a countable set and let $\Lambda$ be a symmetric matrix over $\Ni$ indexed by $I$, with $\Lambda_{i,i}=1$ for all $i\in I$. There exist Cartan masas $(D_i)_{i\in I}$ in the hyperfinite \IIi factor such that $\puk{D_i,D_j}=\{\Lambda_{i,j}\}$ for all $i,j\in I$.
\end{thm}
\begin{proof}
Let $I$ and $\Lambda$ be as in the statement of Theorem \ref{SumCartan.CountableFamily}.  For each unordered pair $\{i,j\}$ of distinct elements of $I$, use Lemma \ref{SumCartan.ABC} to find Cartan masas $(D^{\{i,j\}}_r)_{r\in I}$ in the copy of the hyperfinite \IIi factor denoted $R^{\{i,j\}}$ such that
$$
\puk{D^{\{i,j\}}_r,D^{\{i,j\}}_s}=\begin{cases}\{\Lambda_{i,j}\}&\{r,s\}=\{i,j\}\\\{1\}&\textrm{otherwise}\end{cases}.
$$
This is achieved by taking $D^{(i,j)}_i=A$, $D^{(i,j)}_j=B$ and $D^{(i,j)}_r=C$ for $r\neq i,r\neq j$ where $A,B,C$ are the masas resulting from taking $n=\Lambda_{i,j}$ in Lemma \ref{SumCartan.ABC}.  Now form the copy of the hyperfinite \IIi factor $R=\vnotimes_{\{i,j\}}R^{\{i,j\}}$ and masas $D_r=\vnotimes_{\{i,j\}}D_r^{\{i,j\}}$ for $r\in I$.  Lemma \ref{Prelim.MixTensor} ensures these masas have
$$
\puk{D_i,D_j}=\{\Lambda_{i,j}\}
$$
for all $i,j\in I$.\qed
\end{proof}

We can immediately deduce the existence of masas with certain Puk\'anszky invariants.  The subsets below where first found in \cite{Stormer.Puk} using ergodic methods.
\begin{cor}\label{SumCartan.1InPuk}
Let $E$ be a finite subset of $\Ni$ with $1\in E$.  Then there exists a masa in the hyperfinite \IIi factor whose Puk\'anszky invariant is $E$.
\end{cor}
\begin{proof}
If we work in the $n\times n$ matrices $M_n(R)$ over the hyperfinite \IIi factor, and form the direct sum $A=D_1\oplus D_2\oplus\dots\oplus D_n$ of $n$ Cartan masas, then 
$$
\puk{A}=\{1\}\cup\bigcup_{i<j}\puk{D_i,D_j}.
$$
The corollary then follows from Theorem \ref{SumCartan.CountableFamily} by choosing a large but finite $I$ and appropriate values of $\Lambda_{i,j}$ depending on the set $E$.\qed
\end{proof}

All the pairs of Cartan masas we have produced have had a singleton for their mixed Puk\'anszky invariant.  What are the possible values of $\puk{A,B}$ when $A$ and $B$ are Cartan masas in a \IIi factor? 

\section{The main construction}\label{Values}
In this section we give a construction of masas in McDuff \IIi factors, which we use to establish the main results of the paper in section \ref{Main}.  We need to introduce a not insubstantial amount of notation.  Let $N_0$ be a fixed separable McDuff \IIi factor and for each $r\in\mathbb N$, let $R^{(r)}$ be a copy of the hyperfinite \IIi factor.  Let $N_r=N_0\ \vnotimes\ R^{(1)}\ \vnotimes\ \dots\ \vnotimes\  R^{(r)}$ so that with the inclusion map $x\mapsto x\otimes 1_{R^{(r+1)}}$ we can regard $N_r$ as a von Neumann subalgebra of $N_{r+1}$.  We let $N$ be the direct limit of this chain, so that
$$
N=(N_0\ \vnotimes\ \bigotimes_{r=1}^\infty R^{(r)})''
$$
acting on $L^2(N_0)\otimes\bigotimes_{r=1}^\infty L^2(R^{(r)})$.  The \IIi factor $N$ is isomorphic to $N_0$ and we shall regard all the $N_r$ as subalgebras of $N$.

Whenever we have a masa $D$ inside a \IIi factor, we are able to use the isomorphism between $D$ and $L^\infty[0,1]$ to choose families of projections $e^{(m)}_i(D)$ in $D$ for $m\in\mathbb N$ and $i=(i_1,\dots,i_m)\in\{0,1\}^m$, which satisfy:
\begin{enumerate}
\item For each $m$ the $2^m$ projections $e^{(m)}_i(D)$ are pairwise orthogonal and each projection has trace $2^{-m}$;
\item For each $m$ and $i=(i_1,\dots,i_m)\in\{0,1\}^m$ we have
$$
e_i^{(m)}(D)=e_{i\vee0}^{(m+1)}(D)+e_{i\vee1}^{(m+1)}(D),
$$
where $i\vee 0=(i_1,\dots,i_m,0)$ and $i\vee 1=(i_1,\dots,i_m,1)$;
\item The projections $e^{(m)}_i(D)$ generate $D$.
\end{enumerate}
In the procedure that follows we shall assume that masas come with these projections when needed.

For $m\in\mathbb N$ and $r\geq 0$, let $\ISet{r}{m}$ denote the set of all $i=(i^{(0)},i^{(1)},\dots,i^{(r)})$ where $i^{(r-s)}=(i^{(r-s)}_1,i^{(r-s)}_2,\dots,i^{(r-s)}_{m+s})\in\{0,1\}^{m+s}$ is a sequence of zeros and ones of length $m+s$.  In this way the last sequence, $i^{(r)}$, has length $m$ and each earlier sequence is one element longer than the following sequence.  We have restriction maps from $\ISet{r}{m}$ to $\ISet{r-1}{m+1}$ obtained by forgetting about the last sequence $i^{(r)}$. Note that $i^{(r-1)}$ has length $m+1$ so that this restriction does lie in $\ISet{r-1}{m+1}$.  We can also restrict by shortening the length of all the sequences. In full generality we have restriction maps from $\ISet{r}{m}$ into $\ISet{s}{l}$ whenever $s\leq r$ and $l\leq m+r-s$.  Given $i\in\ISet{r}{m}$ and $k\in\ISet{s}{l}$ (for $s\leq r$ and $l\leq m+r-s$) write $i\geq k$ if the restriction of $i$ to $\ISet{s}{l}$ is precisely $k$.  When $i\in\ISet{r}{m}$ for some $r$, we write $i|_s$ for the restriction of $i$ to $\ISet{s}{1}$ for $s\leq r$.  We take $i|_{-1}=j|_{-1}$ as a convention for all $i,j\in\ISet{r}{m}$.

The inputs to our construction are a masa $A_0$ in $N_0$ and values $\Lambda^{(r)}_{i,j}=\Lambda^{(r)}_{j,i}\in\Ni$ for all $r=0,1,2,\dots$ and $i,j\in\ISet{r}{1}$ with $i\neq j$ and $i|_{r-1}=j|_{r-1}$.   We regard these as fixed henceforth.  For $i\in\ISet{0}{m}$, define $f^{(0,m)}_i=e^{(m)}_{i^{(0)}}(A_0)$.  Suppose inductively that we have produced masas $A_s\subset N_s$ for each $s\leq r$ and that, for each $m\in\mathbb N$, projections $(f^{(s,m)}_i)_{i\in\ISet{s}{m}}$ in $A_s$ have been specified such that:
\begin{enumerate}[(i)]
\item For each $m\in\mathbb N$ and $s\leq r$, the $|\ISet{s}{m}|$ projections $(f^{(s,m)}_i)_{i\in\ISet{s}{m}}$ are pairwise orthogonal and each has trace $|\ISet{m}{s}|^{-1}$;\label{Values.Glue.1}
\item For each $m\in\mathbb N$, $s\leq r$ and $i\in\ISet{s}{m}$ we have
$$
f^{(s,m)}_i=\sum_{\substack{j\in\ISet{s}{m+1}\\j\geq i}}f^{(s,m+1)}_j;\label{Values.Glue.2}
$$
\item For any $s\leq t\leq r$ and $i\in\ISet{s}{m+t-s}$ we have
$$
f^{(s,m+t-s)}_i=\sum_{\substack{j\in\ISet{t}{m}\\j\geq i}}f^{(t,m)}_j,\label{Values.Glue.3}
$$
noting that in this statement we regard the $f^{(s,m+t-s)}$ as lying inside $N_t$;
\item For each $s\leq r$ the projections $\Set{f^{(s,m)}_i|m\in\mathbb N,\quad i\in\ISet{s}{m}}$ generate $A_s$.\label{Values.Glue.4}
\end{enumerate}
Note that conditions (\ref{Values.Glue.3}) and (\ref{Values.Glue.4}) ensure that $A_s\subset A_t$.

To define $A_{r+1}$, use Theorem \ref{SumCartan.CountableFamily} to produce Cartan masas $(D^{(r+1)}_i)_{i\in\ISet{r}{1}}$ in $R^{(r+1)}$ such that when $i\neq j$ we have
\begin{equation}\label{Values.Induct}
\puk{D^{(r+1)}_i,D^{(r+1)}_j}=\begin{cases}\{\Lambda^{(r)}_{i,j}\}&i|_{r-1}=j|_{r-1}\\ \{1\}&\textrm{otherwise}\end{cases}.
\end{equation}
Let $A_{r+1}$ be given by
\begin{equation}\label{Values.Induct.2}
A_{r+1}=\bigoplus_{i\in\ISet{r}{1}} A_rf^{(r,1)}_i\otimes D^{(r+1)}_i
\end{equation}
a masa in $N_r\ \vnotimes\ R^{(r+1)}=N_{r+1}$, which has $A_r\subset A_{r+1}$.  To complete the inductive construction we must define $f^{(r+1,m)}_i$ for $i\in\ISet{r+1}{m}$ in a manner which satisfies conditions (\ref{Values.Glue.1}) through (\ref{Values.Glue.4}) above.  Given $m\in\mathbb N$ and $i\in\ISet{r+1}{m}$, let $i'$ be the restriction of $i$ to $\ISet{r}{m+1}$ and recall that $i|_r$ is the restriction of $i$ to $\ISet{r}{1}$.  Now define
\begin{equation}\label{Values.DefF}
f^{(r+1,m)}_i=f^{(r,m+1)}_{i'}\otimes e^{(m)}_{i^{(r+1)}}(D^{(r+1)}_{i|_r}).
\end{equation}
Since $f^{(r,m+1)}_{i'}\leq f^{(r,1)}_{i|_r}$, this does define a projection in $A_{r+1}$.  That the $f^{(r+1,m)}_i$ satisfy the required conditions is routine.  We give the details as Lemma \ref{Values.Glue} below for completeness.
\begin{lem}\label{Values.Glue}
The projections $(f^{(r+1,m)}_i)_{i\in\ISet{r+1}{m}}$ defined in (\ref{Values.DefF}) satisfy the conditions (\ref{Values.Glue.1}) through (\ref{Values.Glue.4}) above.
\end{lem}
\begin{proof}
For $m\in\mathbb N$ fixed, the projections $(f^{(r+1,m)}_i)_{i\in\ISet{r+1}{m}}$ are pairwise orthogonal and have trace $|\ISet{r+1}{m}|^{-1}$ as the projections $(f^{(r,m+1)}_{i'})_{i'\in\ISet{r}{m+1}}$ are pairwise orthogonal with trace $|\ISet{r}{m+1}|^{-1}$ and the projections $(e^{(m)}_{j}(D^{(r+1)}_{i|_r}))_{j\in\{0,1\}^m}$ are also pairwise orthogonal and each have trace $2^{-m}$.  In this way the projections for $A_{r+1}$ satisfy condition (\ref{Values.Glue.1}).

For condition (\ref{Values.Glue.2}), fix $i\in\ISet{r+1}{m}$ for some $m\in\mathbb N$ and let $i'$ be as in the definition of $f^{(r+1,m)}_i$.  Now
\begin{align*}
f^{(r+1,m)}_i=&f^{(r,m+1)}_{i'}\otimes e^{(m)}_{i^{(r+1)}}(D^{(r+1)}_{i|_r})\\
=&\sum_{\substack{j'\in\ISet{r}{m+2}\\j'\geq i'}}f^{(r,m+2)}_{j'}\otimes \Bigg(e^{(m+1)}_{i^{(r+1)}\vee 0}(D^{(r+1)}_{i|_r})+e^{(m+1)}_{i^{(r+1)}\vee 1}(D^{(r+1)}_{i|_r})\Bigg)\\
=&\sum_{\substack{j\in\ISet{r+1}{m+1}\\j\geq i}}f^{(r+1,m+1)}_j
\end{align*}
from condition (\ref{Values.Glue.2}) for the $f^{(r,m+1)}_{i'}$ and the second condition in the definition of the $e^{(m)}_{k}(D)$.  This is precisely condition (\ref{Values.Glue.2}).

We only need to check condition (\ref{Values.Glue.3}) when $t=r+1$, so take $s\leq r$, $m\in\mathbb N$ and $i\in\ISet{s}{m+r+1-s}$.  By the inductive version of (\ref{Values.Glue.3}) we have
$$
f^{(s,m+r+1-s)}_i=\sum_{\substack{j\in\ISet{r}{m+1}\\j\geq i}}f^{(r,m+1)}_j.
$$
For each $j\in\ISet{r}{m+1}$ with $j\geq i$ we have
\begin{align*}
f^{(r,m+1)}_j\otimes 1_{R^{(r+1)}}=&f^{(r,m+1)}_j\otimes\sum_{j^{(r+1)}\in\{0,1\}^m}e^{(m)}_{j^{(r+1)}}(D^{(r+1)}_{j|_r})\\
=&\sum_{\substack{k\in\ISet{r+1}{m}\\k\geq j}}f^{(r+1,m+1)}_k,
\end{align*}
where $j|_r$ is the restriction of $j$ to $\ISet{r}{1}$.  Therefore,
$$
f^{(s,m+r+1-s)}_i=\sum_{\substack{k\in\ISet{r+1}{m}\\k\geq i}}f^{(r+1,m+1)}_k,
$$
which is condition (\ref{Values.Glue.3}).

For $j\in\ISet{r}{1}$, the projections $f^{(r,m)}_k$ indexed by $k\in\ISet{r}{m}$ with $k\geq j$ generate the cut-down $A_rf^{(r,1)}_j$.  Hence the projections $f^{(r+1,m)}_i$, for $i\in\ISet{r+1}{m}$ with $i\geq j$ generate $A_rf^{(r,1)}_j\otimes D^{(r+1)}_j$.  In this way we see that the projections $f^{(r+1,m)}_i$ for $i\in\ISet{r+1}{m}$ generate $A_{r+1}$, which is condition (\ref{Values.Glue.4}).\qed
\end{proof}

This completes the inductive stage of the construction.  We have masas $A_r$ in $N_r$ for each $r$ such that $A_r\otimes 1_{R^{(r+1)}}\subset A_{r+1}$.  We shall regard all these masas as subalgebras of the infinite tensor product \IIi factor $N$, where they are no longer maximal abelian.  Define $A=(\bigcup_{r=0}^\infty A_r)''$, an abelian subalgebra of $R$.  For $r\geq 0$ we have
$$
A_r'\cap N=A_r\ \vnotimes\ R^{(r+1)}\ \vnotimes\ R^{(r+2)}\ \vnotimes\ \dots
$$
so that for $x\in N_r\subset N$ we have $\ce{A_r'\cap N}{x}=\ce{A_r}{x}$, where $\mathbb E_M$ denotes the unique trace-preserving conditional expectation onto the von Neumann subalgebra $M$.  As $A_r\subset A\subset A'\cap N\subset A_r'\cap N$ we obtain $\ce{A}{x}=\ce{A'\cap N}{x}$ for any $x\in\bigcup_{r=0}^\infty N_r$.  These $x$ are weakly dense in $N$ so $A=A'\cap N$ is a masa in $N$, see \cite[Lemma 2.1]{Popa.Kadison}.

 \section{The Puk\'anszky invariant of $A$}\label{MainLem}
Our objective here is to compute the Puk\'anskzy invariant of the masas of section \ref{Values} in terms of the masa $A_0$ and the specified values $\Lambda_{i,j}^{(r)}$.  Following the usual convention, we shall write $\A$ for the algebra $(A\cup JAJ)''$, an abelian subalgebra of $\mathbb B(L^2(N))$.
\begin{lem}\label{Values.Diagonal}
Let $A$ be a masa produced by means of the construction described in section \ref{Values}.  Then 
$$
\puk{A}=\bigcup_{r=0}^\infty\bigcup_{\substack{i,j\in\ISet{r}{1}\\i\neq j\\i|_{r-1}=j|_{r-1}}}\type{\A'f^{(r,1)}_iJf^{(r,1)}_jJ}.
$$
\end{lem}
\begin{proof}
Fix $s\geq 0,m\in\mathbb N$ and $i\in\ISet{s}{m}$.  Let $r=s+m-1$, so that condition (\ref{Values.Glue.3}) gives
$$
f^{(s,m)}_i=\sum_{\substack{j\in\ISet{r}{1}\\j\geq i}}f^{(r,1)}_j.
$$
Condition (\ref{Values.Glue.4}) shows that the projections $f^{(s,m)}_i$, for $m\in\mathbb N$ and $i\in\ISet{s}{m}$, generate $A_s$.  Hence every $A_s$ is contained in the abelian von Neumann algebra generated by all the $f^{(r,1)}_i$ for $i\in\ISet{r}{1}$ and $r\geq 0$, so these projections generate $A=(\bigcup_{s=1}^\infty A_s)''$.

Writing $B_r$ for the abelian von Neumann subalgebra of $N$ generated by the projections $(f^{(r,1)}_i)_{i\in\ISet{r}{1}}$, Lemma 2.1 of \cite{Popa.Kadison} shows us that 
$$
\lim_{r\rightarrow\infty}\nm{\ce{B_r'\cap N}{x}-\ce{A}{x}}_2=0
$$
for all $x\in N$, where $\mathbb E_M$ denotes the trace-preserving conditional expectation onto the von Neumann subalgebra $M$ of $N$.  It is well known that $\mathbb E_{B_r'\cap N}=\sum_{i\in\ISet{r}{1}}f^{(r,1)}_iJf^{(r,1)}_iJ$ in this case, so
$$
e_A=\lim_{r\rightarrow\infty}\sum_{i\in\ISet{r}{1}}f^{(r,1)}_iJf^{(r,1)}_iJ,
$$
with strong-operator convergence.  Hence
$$
1-e_A=\sum_{r=0}^\infty\sum_{\substack{i,j\in\ISet{r}{1}\\i\neq j\\i|_{r-1}=j|_{r-1}}}f^{(r,1)}_iJf^{(r,1)}_jJ
$$
so the only contributions to the Puk\'anskzy invariant of $A$ come from the central cutdowns $\A'f^{(r,1)}_iJf^{(r,1)}_jJ$ for $r\geq 0$, $i,j\in\ISet{r}{1}$ with $i\neq j$ and $i|_{r-1}=j|_{r-1}$.\qed
\end{proof}

For $s\geq 0$, write $\A_s$ for the abelian von Neumann algebra $(A_s\cup JA_sJ)''\subset\mathbb B(L^2(N_s))$.  For the rest of this section we shall denote operators in $\mathbb B(L^2(N_s))$ with a superscript $^{(s)}$.  Since 
$$
\mathbb B(L^2(N_{s+1}))=\mathbb B(L^2(N_s))\ \vnotimes\ \mathbb B(L^2(R^{(s+1)}))
$$
we have $T^{(s)}\otimes I_{L^2(R^{(s+1)})}\in\mathbb B(L^2(N_{s+1}))$ for all $T^{(s)}\in\mathbb B(L^2(N_s))$.  We shall write $T^{(s+1)}$ for this operator, and 
$$
T=T^{(s)}\otimes I_{L^2(R^{(s+1)})}\otimes I_{L^2(R^{(s+2)})}\otimes\dots
$$
for this extension of $T^{(s)}$ to $L^2(N)$. We refer to these operators as the canonical extensions of $T^{(s)}$.  For $T^{(s)}\in\A_s$, we have $T^{(s+1)}\in\A_{s+1}$ and $T\in\A$, since $A_s\subset A_{s+1}\subset A$.  Let $p_s$ denote the orthogonal projection from $L^2(N)$ onto $L^2(N_s)$.
\begin{prop}\label{MainLem.Easy}
Let $s\geq 0$ and $T^{(s)}\in\mathbb B(L^2(N_s))$.  Then $T^{(s)}\in\A_s'$ if and only if the extension $T$ lies in $\A'$.  Also $p_s\A'p_s=\A_s'$.
\end{prop}
\begin{proof}
Let $T\in\mathbb B(L^2(N))$ lie in $\A'$.  For each $s$ and $x\in A_s$, we have $p_sxp_s=xp_s=p_sx$ and $p_sJxJp_s=JxJp_s=p_sJxJ$.  Then $p_sTp_s$ commutes with both $x$ and $JxJ$ and hence lies in $\A_s'$.  This gives $p_s\A'p_s\subset\A_s'$ and shows that if $T$ is the canonical extension of some $T^{(s)}\in\mathbb B(L^2(N_s))$, then $T^{(s)}\in\A_s'$.

For the converse, consider $T^{(s)}\in\A_s'$ and take $x\in A_{s+1}$ so that
$$
x=\sum_{i\in\ISet{s}{1}}x_if^{(s,1)}_i\otimes y_i
$$
for some $x_i\in A_s$ and $y_i\in D^{(s+1)}_i$ by the inductive definition of $A_{s+1}$ in equation (\ref{Values.Induct.2}).  Then $T^{(s+1)}$ commutes with $x$ since $T^{(s)}$ commutes with each $x_if^{(s,1)}_i$.  Similarly $T^{(s+1)}$ commutes with $JxJ$, so $T^{(s+1)}\in\A_{s+1}'$.  Proceeding by induction, we see that $T^{(r)}\in\A_r'$ for all $r\geq s$.  Hence, the canonical extension $T$ commutes with $x$ and $JxJ$ for all $x\in\bigcup_{r=0}^\infty A_r$ and these elements are weakly dense in $\A$, so $T\in\A'$. For $T^{(s)}\in\mathbb B(L^2(N_s))$ the canonical extension $T$ has $p_sTp_s=T^{(s)}$, so $\A_s'\subset p_s\A'p_s$.\qed
\end{proof}

Our objective is to determine the type decomposition of the $\A'f^{(r,1)}_iJf^{(r,1)}_jJ$ appearing in Lemma \ref{Values.Diagonal}.  For $r\geq 0$ and $i\in\ISet{r}{1}$, the inductive definition (\ref{Values.DefF}) ensures that
$$
f^{(r,1)}_i=e^{(r+1)}_{i^{(0)}}(A_0)\otimes e^{(r)}_{i^{(1)}}(D^{(1)}_{i|_0})\otimes\dots\otimes e^{(1)}_{i^{(r)}}(D^{(r)}_{i|_{r-1}})
$$
recalling that $i|_s$ is the restriction of $i$ to $\ISet{s}{1}$.
\begin{lem}\label{MainLem.Full}
Let $r\geq 0$ and $i,j\in\ISet{r}{1}$ have $i\neq j$ and $i|_{r-1}=j|_{r-1}$. Let $Q^{(0)}\in\A_0e^{(r+1)}_{i^{(0)}}(A_0)Je^{(r+1)}_{j^{(0)}}(A_0)J$ be a non-zero projection such that $\A_0'Q^{(0)}$ is homogeneous of type $\mathrm{I}_m$ for some $m\in\Ni$.  Then, writing $Q$ for the canonical extension of $Q^{(0)}$ to $L^2(N)$, $\A'f^{(r,1)}_iJf^{(r,1)}_jJQ$ is homogeneous of type $\mathrm{I}_{m\Lambda^{(r)}_{i,j}}$.
\end{lem}
\begin{proof}
Fix $m\in\Ni$ and $Q^{(0)}\neq 0$ as in the statement of the Lemma. Observe that
\begin{align*}
A_{r+1}f^{(r,1)}_i=&A^{(r)}f^{(r,1)}_i\ \vnotimes\  D^{(r+1)}_i\\
=&A_0e^{(r+1)}_{i^{(0)}}(A_0)\ \vnotimes\  D^{(1)}_{i|_0}e^{(r)}_{i^{(1)}}(D^{(1)}_{i|_0})\ \vnotimes\ \dots\ \vnotimes\  D^{(r)}_{i|_{r-1}}e^{(1)}_{i^{(r)}}(D^{(r)}_{i|_{r-1}})\ \vnotimes\ D^{(r+1)}_i
\end{align*}
so that
\begin{align*}
&\A_{r+1}f^{(r,1)}_iJf^{(r,1)}_jJQ^{(r+1)}\\
=&\A_0Q^{(0)}\ \vnotimes\ (D^{(1)}_{i|_0}\cup JD^{(1)}_{i|_0}J)''e^{(r)}_{i^{(1)}}(D^{(1)}_{i^{(0)}})Je^{(r)}_{j^{(1)}}(D^{(1)}_{i^{(0)}})J\\
&\quad\vnotimes\ \dots\ \vnotimes\ (D^{(r)}_{i|_{r-1}}\cup JD^{(r)}_{i|_{r-1}}J)''e^{(1)}_{i^{(r)}}(D^{(r)}_{i|_{r-1}})Je^{(1)}_{j^{(r)}}(D^{(r)}_{i|_{r-1}})J\ \vnotimes\ (D^{(r+1)}_i\cup JD^{(r+1)}_jJ)'',
\end{align*}
using $i|_s=j|_s$ for $s=0,\dots,r-1$.  We are also abusing notation by writing $J$ for the modular conjugation operator regardless of the space on which it operates.  Taking commutants gives
\begin{align*}
&\A_{r+1}'f^{(r,1)}_iJf^{(r,1)}_jJQ^{(r+1)}\\
=&\A_0'Q^{(0)}\ \vnotimes\ (D^{(1)}_{i|_0}\cup JD^{(1)}_{i|_0}J)'e^{(r)}_{i^{(1)}}(D^{(1)}_{i^{(0)}})Je^{(r)}_{j^{(1)}}(D^{(1)}_{i^{(0)}})J\\
&\quad\vnotimes\ \dots\ \vnotimes\ (D^{(r)}_{i|_{r-1}}\cup JD^{(r)}_{i|_{r-1}}J)'e^{(1)}_{i^{(r)}}(D^{(r)}_{i|_{r-1}})Je^{(1)}_{j^{(r)}}(D^{(r)}_{i|_{r-1}})J\ \vnotimes\ (D^{(r+1)}_i\cup JD^{(r+1)}_jJ)'.
\end{align*}
For $s\leq r$, each $(D^{(s)}_{i|_{s-1}}\cup JD^{(s)}_{i|_{s-1}}J)''$ is maximal abelian in $\mathbb B(L^2(R^{(s)}))$ since $D^{(s)}_{i|_{s-1}}$ is a Cartan masa so has Puk\'anszky invariant $\{1\}$.  The masas $D^{(r+1)}_k$ where defined in (\ref{Values.Induct}) so that $(D^{(r+1)}_i\cup JD^{(r+1)}_jJ)'$ is homogeneous of type $\mathrm{I}_{\Lambda^{(r)}_{i,j}}$.  We learn that $\A_{r+1}'f^{(r,1)}_iJf^{(r,1)}_jJQ^{(r+1)}$ is homogeneous of type $\mathrm{I}_{m'}$, where $m'=m\Lambda_{i,j}^{(r)}$.

Find a family of pairwise orthogonal projections $(Q^{(r+1)}_q)_{0\leq q<m'}$ with sum $Q^{(r+1)}$ and which are equivalent abelian projections in $\A_{r+1}'f^{(r,1)}_iJf^{(r,1)}_jJQ^{(r+1)}$.  The canonical extensions $(Q_q)_{0\leq q<m'}$ to $L^2(N)$ form a family of pairwise orthogonal projections in $\A'Q$ (by Proposition \ref{MainLem.Easy}) with sum $Q$.  These projections are equivalent in $\A'Q$ as if $V^{(r+1)}$ is a partial isometry in $\A_r'Q^{(r+1)}$ with $V^{(r+1)}{V^{(r+1)}}^*=Q_q$ and ${V^{(r+1)}}^*V^{(r+1)}=Q_{q'}$, then Proposition \ref{MainLem.Easy} ensures that the canonical extension $V$ lies in $\A'$.  It is immediate that $VV^*=Q_q$ and $V^*V=Q_{q'}$.  We shall show that these projections are abelian projections in $\A'$.  It will then follow that $\A'Q$ is homogeneous of type $\mathrm{I}_{m'}$.

For $s\geq r+1$ and $k,l\in\ISet{s}{1}$ with $k\geq i$ and $l\geq j$, we have 
$$
A_{s+1}f^{(s,1)}_k=A_sf^{(s,1)}_k\ \vnotimes\ D^{(s+1)}_k
$$
so that
$$
\A_{s+1}(f^{(s,1)}_kJf^{(s,1)}_lJ)Q^{(s+1)}\cong \A_s(f^{(s,1)}_kJf^{(s,1)}_lJ)Q^{(s)} \vnotimes\ (D_k^{(s+1)}\cup JD_l^{(s+1)}J)''.
$$
Again we take commutants to obtain
$$
\A_{s+1}'(f^{(s,1)}_kJf^{(s,1)}_lJ)Q^{(s+1)}\cong \A_s'(f^{(s,1)}_kJf^{(s,1)}_lJ)Q^{(s)}\ \vnotimes\  (D_k^{(s+1)}\cup JD_l^{(s+1)}J)'.
$$
Since $i\neq j$ it is not possible for $k|_s$ to equal $l|_s$, so (\ref{Values.Induct}) shows us that $(D_k^{(s+1)}\cup JD_l^{(s+1)}J)'$ is abelian.  Therefore, if $Q^{(s)}_qf^{(s,1)}_kJf^{(s,1)}_lJ$ (some $q=1,\dots,m'$) is an abelian projection in $\A_s'$, then $Q^{(s+1)}_qf^{(s+1,1)}_kJf^{(s+1,1)}_lJ$ is abelian in $\A_{s+1}'$.  The projections $f^{(s,1)}_kJf^{(s,1)}_lJ$ are central and satisfy
$$
\sum_{\substack{k,l\in\ISet{s}{1}\\k\geq i\\l\geq j}}f^{(s,1)}_kJf^{(s,1)}_lJ=f^{(r,1)}_iJf^{(r,1)}_jJ.
$$
By induction and summing over all $k\geq i$ and $l\geq j$, we learn that $(Q_q^{(s)})_{0\leq q<m'}$ form a family of equivalent abelian projections in $\A'Q^{(s)}$ with sum $s$ for every $s\geq r+1$.  

For $s\geq r+1$ and each $q$, the algebras $\A_s'Q^{(s)}_q=p_s\A'Q_qp_s$ are abelian.  Since the projections $p_s$ tend strongly to the identity, we see that each $\A'Q_q$ is abelian too.\qed
\end{proof}

We can now describe the Puk\'anszky invariant of the masas in section \ref{Values}. 
\begin{thm}\label{MainLem.Puk}
Let $A$ be a masa in a separable McDuff \IIi factor produced via the construction of section \ref{Values}.  That is we are given a masa $A_0\subset N_0$ and values $\Lambda_{i,j}^{(r)}=\Lambda_{j,i}^{(r)}\in\Ni$ for $r\geq 0$, $i,j\in\ISet{r}{1}$ with $i\neq j$ and $i|_{r-1}=j|_{r-1}$.  Then
\begin{equation}\label{MainLem.Puk.1}
\puk{A}=\bigcup_{r=0}^\infty\bigcup_{\substack{i,j\in\ISet{r}{1}\\i\neq j\\i|_{r-1}=j|_{r-1}}}\Lambda^{(r)}_{i,j}\cdot\type{\A_0'e^{(r+1)}_{i^{(0)}}(A_0)Je^{(r+1)}_{j^{(0)}}(A_0)J}.
\end{equation}
\end{thm}
\begin{proof}
For $r\geq 0$, $i,j\in\ISet{r}{1}$ with $i\neq j$ and $i|_{r-1}=j|_{r-1}$, it follows from Lemma \ref{MainLem.Full} that
$$
\type{\A'f^{(r,1)}_iJf^{(r,1)}_jJ}=\Lambda^{(r)}_{i,j}\cdot\type{\A_0'e^{(r+1)}_{i^{(0)}}(A_0)\cup Je^{(r+1)}_{j^{(0)}}(A_0)J}.
$$
The theorem then follows from Lemma \ref{Values.Diagonal}.\qed
\end{proof}

\section{Main results}\label{Main}
We start by applying Theorem \ref{MainLem.Puk} when $\puk{A_0}$ is a singleton.
\begin{thm}
For $n\in\mathbb N$, suppose that $N_0$ is a separable McDuff \IIi factor containing a masa with Puk\'anszky invariant $\{n\}$.  For every non-empty set $E\subset\Ni$, there exists a masa $A$ in $N_0$ with $\puk{A}=\{n\}\cdot E$.
\end{thm}
\begin{proof}
Let $A_0$ be a masa in $N_0$ with $\puk{A}=\{n\}$ and choose the values $\Lambda_{i,j}^{(r)}=\Lambda_{j,i}^{(r)}$ for $r\geq 0$ and $i,j\in\ISet{r}{1}$ with $i\neq j$ and $i|_{r-1}=j|_{r-1}$ so that
$$
E=\Set{\Lambda^{(r)}_{i,j}|r\geq 0,\quad i,j\in\ISet{r}{1},\quad i\neq j,\quad i|_{r-1}=j|_{r-1}}.
$$
The resulting masa $A$ in $N\cong N_0$ produced by the main construction has Puk\'anszky invariant $\{n\}\cdot E$ by Theorem \ref{MainLem.Puk}.\qed
\end{proof}

Since Cartan masas have Puk\'anskzy invariant $\{1\}$, we obtain the following Corollary immediately.  
\begin{cor}\label{Main.SimpleMcDuff}
Let $N$ be a McDuff \IIi factor containing a simple masa, for example a Cartan masa.  Every non-empty subset of $\Ni$ arises as the Puk\'anszky invariant of a masa in $N$.
\end{cor}

A little more care enables us to address the question of the range of the Puk\'anszky invariant on singular masas in the hyperfinite \IIi factor and other McDuff \IIi factors containing a simple singular masa.  Puk\'anszky's original work \cite{Pukanszky.Invariant} exhibits a simple singular masa in the hyperfinite \IIi factor.
\begin{cor}\label{Main.Singular}
Let $N$ be a separable McDuff factor containing a simple singular masa, such as the hyperfinite \IIi factor.  Given any non-empty $E\subset\Ni$ there is a singular masa $A$ in $N$ with $\puk{A}=E$.
\end{cor}
\begin{proof}
If $1\not\in E$, a masa in $N$ with Puk\'anszky invariant $E$ is automatically singular by \cite[Remark 3.4]{Popa.NotesCartan}.  We have already produced these masas in Corollary \ref{Main.SimpleMcDuff}.  The hypothesis ensures us a simple singular masa in $N$.  For the remaining case of some $E\neq \{1\}$ with $1\in E$, let $A_1$ be a singular masa in $N$ with $\pukN{A_1}{N_1}=\{1\}$ and $A_2$ be a singular masa in the hyperfinite \IIi factor $R$ with $\pukN{A_2}{R}=E\setminus\{1\}$.  Then $A=A_1\vnotimes A_2$ is a masa in $N\vnotimes R\cong N$. Lemma 2.1 of \cite{Sinclair.Puk} ensures that
$$
\puk{A}=\{1\}\cup (E\setminus \{1\})\cup 1\cdot(E\setminus\{1\})=E.
$$
The singularity of $A$ is Corollary 2.4 of \cite{Saw.StrongSing}.\qed
\end{proof}

Next we justify the claims made at the end of section \ref{Prelim}.
\begin{thm}\label{Main.DirectSum}
Let $E,F,G\subset\Ni$ be non-empty.  Then there exist masas $B$ and $C$ in the hyperfinite \IIi factor with $\puk{B}=E$, $\puk{C}=F$ and $\puk{B,C}=G$.
\end{thm}
\begin{proof}
Let $R_0$ be a copy of the hyperfinite \IIi factor and $A_0$ a Cartan masa in $R_0$.  An element $k$ of $\ISet{0}{1}$ is of the form $(k^{(0)})$ where $k^{(0)}$ is a $1$-tuple --- either $0$ or $1$.  Write $\mathbf{0}$ and $\mathbf{1}$ for these two elements and let $e_0=f^{(1)}_\mathbf{0}$ and $e_1=f^{(1)}_\mathbf{1}$ so that $e_0$ and $e_1$ are orthogonal projections in $A$ with $\tr(e_0)=\tr(e_1)=1/2$.  Choose the $\Lambda^{(r)}_{i,j}=\Lambda_{j,i}^{(r)}$ such that:
\begin{align*}
E=&\Set{\Lambda^{(r)}_{i,j}|r\geq 1,\quad i,j\in\ISet{r}{1},i\neq j,i|_{r-1}=j|_{r-1},i,j\geq \mathbf{0}},\\
F=&\Set{\Lambda^{(r)}_{i,j}|r\geq 1,\quad i,j\in\ISet{r}{1},i\neq j,i|_{r-1}=j|_{r-1},i,j\geq \mathbf{1}},\\
G=&\Set{\Lambda^{(r)}_{i,j}|r\geq 0,\quad i,j\in\ISet{r}{1},i\neq j,i|_{r-1}=j|_{r-1}, i\geq \mathbf{0},j\geq \mathbf{1}}\\
=&\Set{\Lambda^{(r)}_{i,j}|r\geq 0,\quad i,j\in\ISet{r}{1},i\neq j,i|_{r-1}=j|_{r-1}, i\geq\mathbf{1},j\geq\mathbf{0}}.
\end{align*}
For $r,s=0,1$, let $Q_{r,s}=(1-e_A)e_rJe_sJ$ a projection in $\A$.  Now Lemma \ref{MainLem.Full} and Lemma \ref{Values.Diagonal} ensure that $\A'Q_{0,0}$ has a non-zero $\mathrm{I}_m$ cutdown if and only if $m\in E$, $\A'Q_{1,1}$ has a non-zero $\mathrm{I}_m$ cutdown if and only if $m\in F$, $\A'(Q_{0,1}+Q_{1,0})$ has a non-zero $\mathrm{I}_m$ cutdown if and only if $m\in G$.

We now regard $A$ as a direct sum.  Consider the copy of the hyperfinite \IIi factor $S=e_0Re_0$ so that choosing a partial isometry $v\in R$ with $v^*v=e_0$ and $vv^*=e_1$ gives rise to an isomorphism between $R$ and $M_2(S)$ --- the $2\times 2$ matrices over $S$.  Define masas in $S$ by $B=Ae_0$ and $C=v^*(Ae_1)v$.  The discussion above ensures that $\puk{B}=E$, $\puk{C}=F$ and $\puk{B,C}=G$.  Note that $\puk{B,C}$ is independent of $v$ by Proposition \ref{Prelim.Inner}.\qed
\end{proof}

\begin{rem}
If $E\subset\Ni$ contains at least two elements then we can modify the construction in section \ref{Values} to produce uncountably many pairwise non-conjugate masas in the hyperfinite \IIi factor $R$ each with Puk\'anszky invariant $E$.  The idea is to control the supremum of the trace of a projection in the masa $A$ such that $\pukN{Ae}{eRe}=\{n\}$ for some fixed $n\in E$.  For each $t\in(0,1)$, we can produce masas $A$ in $R$ and a projection $e\in A$ with $\tr(e)=t$ such that (with the intuitive diagrams of the introduction) the multiplicity structure of $\A$ is represented by Figure \ref{Main.Fig}, with $1$ down the diagonal and $E\setminus\{n\}$ occurring in the unmarked areas.  All these masas must be pairwise non-conjugate.
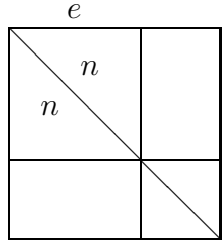
\begin{figure}[h!]
\begin{center}
\begin{picture}(80,90)
\put(0,0){\line(1,0){80}}
\put(0,30){\line(1,0){80}}
\put(0,80){\line(1,0){80}}
\put(0,0){\line(0,1){80}}
\put(50,0){\line(0,1){80}}
\put(80,0){\line(0,1){80}}
\put(0,80){\line(1,-1){80}}
\put(12,47){$n$}
\put(27,62){$n$}
\put(22,84){$e$}
\end{picture}
\caption{The multiplicity structure of $\A$.}
\label{Main.Fig}
\end{center}
\end{figure}

No modifications are required to obtain any diadic rational for $t$, we follow Theorem \ref{Main.DirectSum} to control the multiplicity structure of $\A$.  For general $t$ we can approximate the required structure using diadic rationals, leaving the area we are unable to handle at each stage with multiplicity $1$ so it can be adjusted at a subsequent stage.
\end{rem}

\begin{rem}
For a masa $A$ in a property $\Gamma$-factor $N$, the property that $A$ contains non-trivial centralising sequences for $N$ has been used to distinguish between non-conjugate masas, see for example \cite{Jones.PropertiesMasas,Stormer.Puk,Saw.Path}.  We can easily adjust the construction of section \ref{Values} to ensure that all the masas produced have this property. Suppose that we identify each $R^{(r)}$ with $R^{(r)}\vnotimes R^{(r)}$ and we replace the masas $D^{(r)}_i$ in $R^{(r)}$ by $D^{(r)}_i\vnotimes E^{(r)}$ where $E^{(r)}$ is a fixed Cartan masa in $R^{(r)}$.  By Lemma \ref{Prelim.MixTensor} this does not alter the mixed Puk\'anszky invariants of the family, so the Puk\'anszky invariant of the masa resulting from the construction remains unchanged.  This masa now contains non-trivial centralising sequences for $N$.  By way of contrast, the examples in \cite{Sinclair.Puk,Sinclair.FreePuk} arise from inclusions $H\subset G$ of a an abelian group inside a discrete I.C.C. group $G$ with $gHg^{-1}\cap H=\{1\}$ for all $g\in G\setminus H$.  The resulting masa $\VN{H}$ can not contain non-trivial centralising sequences for the \IIi factor $\VN{G}$, \cite{Popa.Orth}.
\end{rem}

Very recently Ozawa and Popa have shown that not every McDuff \IIi factor contains a Cartan masa. Indeed in \cite{Ozawa.OneCartan} they show that there are no Cartan masas in $\mathcal L{\mathbb F_2}\vnotimes R$. It is not known whether every McDuff factor must contain a simple masa (one with Puk\'anskzy invariant $\{1\}$) or a masa whose Puk\'anszky invariant is a finite subset of $\mathbb N$.  We can however obtain subsets containing $\infty$ as Puk\'anszky invariants of masas in a general separable McDuff \IIi factor.
\begin{thm}\label{Main.GeneralMcDuff}
Let $N$ be a separable McDuff \IIi factor.  For every set $E\subset\Ni$ with $\infty\in E$ there is a singular masa $B$ in $N$ with $\puk{B}=E$.
\end{thm}
\begin{proof}
Taking all the $\Lambda_{i,j}^{(r)}=\infty$, gives us a masa $A$ in $N$ with $\puk{A}=\{\infty\}$ by Theorem \ref{MainLem.Puk} (regardless of the masa $A_0$).  Now use the isomorphism $N\cong N\vnotimes R$, where $R$ is the hyperfinite \IIi factor.  Let $B=A\vnotimes A_1$, where $A_1$ is a singular masa in $R$ with $\pukN{A_1}{R}=E$.  Lemma 2.1 of \cite{Sinclair.Puk} gives
$$
\puk{B}=\{\infty\}\cup E\cup \{\infty\}\cdot E=E.\qed
$$
\end{proof}
In particular every separable McDuff \IIi factor contains uncountably many pairwise non-conjugate singular masas.

\bigskip
\noindent Stuart White\\
Department of Mathematics,\\
University of Glasgow,\\
University Gardens,\\
Glasgow.\\
G12 8QW\\
U.K.\\
\texttt{s.white@maths.gla.ac.uk}

\end{document}